\documentclass[12pt]{amsart}
\usepackage{amsthm}
\newtheorem{theorem}{Theorem}
\newtheorem{corollary}{Corollary}

\newtheorem*{conjecture}{Conjecture}

\theoremstyle{remark}

\begin{document}
\def\Rea{\mathrm{Re}\, }
\def\Ima{\mathrm{Im}\, }
\def\C{\mathbf{C}}
\def\bC{\mathbf{\overline{C}}}
\def\R{\mathbf{R}}
\def\Z{\mathbf{Z}}
\title[Metrics with 4 conic singularities]{Metrics of constant
positive curvature with four conic singularities
on the sphere}
\author{Alexandre Eremenko}\thanks{Supported by NSF grant DMS-1665115.}
\begin{abstract}
We show that for given four points on the sphere and prescribed
angles
at these points which are not multiples of $2\pi$,
the number of metrics of curvature $1$ having conic
singularities with these angles at these points
is finite.


\end{abstract}

\maketitle

\section{Introduction}

We consider metrics of constant positive curvature 
with $n$ conic singularities on the sphere. Without loss of generality
we assume that curvature equals $1$. Such a metric can be described
by the length element $\rho(z)|dz|$, where $z$ is a local conformal
coordinate and $\rho$ satisfies
$$\Delta\log\rho+\rho^2=2\pi\sum_{j=0}^{n-1}(\alpha_j-1)\delta_{a_j},$$
where $a_j$ are the singularities with angles $2\pi\alpha_j$.
The problem is to describe the set of such metrics with prescribed
singularities and angles.

For the recent results
on the problem we refer to \cite{E2}, \cite{ET}, \cite{LW}, \cite{M},
\cite{MP1}, \cite{MP2}.
It is believed that when none of the $\alpha_j$ is an integer, 
the number of such metrics with
prescribed $a_j$ and $\alpha_j$ is finite.
This number has been found in some special cases,
\cite{T}, \cite{E}, \cite{EGT},
\cite{ET}, \cite{LW}, \cite{LT},
in particular, there is at most one such metric in the following two
cases: a) when $n\leq 3$ and the angles are not multiples of $2\pi$,
\cite{T}, \cite{E}, and b)
when $\alpha_j<1$ for all $j$, and $n$ is arbitrary, 
\cite{LT}. A new approach to the uniqueness was recently
published in \cite{Ba1} where a different proof of the main result of \cite{LT}
is obtained. However for large angles and $n\geq 4$
there is usually more than one metric
\cite{CL,EGT,ET}.

In this paper we address the case $n=4$.
We briefly recall the reduction of the problem to a problem
about Heun's equation, see \cite{EGT}.

If $S$ is the sphere equipped with such a metric, one can consider
a developing map $f:S\to\bC$, where $\bC$ is the Riemann sphere equipped
with the standard metric of curvature $1$. Strictly speaking, $f$ is
defined on the universal cover of $S\backslash\{ a_0,\ldots,a_{n-1}\},$ but we
prefer to consider $f$ as a multi-valued function 
with branching at the singularities.
One can write $f=w_1/w_2$ where $w_1$ and $w_2$ are two
linearly independent solutions of the Heun equation, a Fuchsian
equation with four singularities. These four singularities
are the singularities
of the metric, and the angles at the singularities are $2\pi$ times
the exponent differences. Heun's equation can be written as
\begin{equation}
\label{heun0}
w''+\left(\sum_{j=0}^2\frac{1-\alpha_j}{z-a_j}\right)w'+
\frac{Az-q}{(z-a_0)(z-a_1)(z-a_2)}w=0,
\end{equation}
where the singularities are $a_0,a_1,a_2,\infty$, the angles
are $2\pi\alpha_j,\; 0\leq j\leq 3$,
and
$$A=(2+\alpha_3-\alpha_0-\alpha_1-\alpha_2)(2-\alpha_3-\alpha_0-\alpha_1-\alpha_2)/4.$$
Three singularities can be placed at arbitrary points, so
one can choose, for example $(a_0,a_1,a_2)=(0,1,t)$. So for given singularities
and angles,
the set of Heun's equations essentially depends on $2$ parameters: $t$
which describes the quadruple of singularities up to conformal equivalence,
and 
$q$ which is called the {\em accessory parameter}. 

The metric and the differential equation (\ref{heun0})
can be lifted
to a torus via the two-sheeted covering ramified over the four singular
points.
Assuming that the singularities are at $e_1,e_2,e_3,\infty,$
where $e_1+e_2+e_3=0$, we consider the Weierstrass
function 
$$\wp:\C\to S,$$
with primitive periods 
$\omega_1,\omega_2.$ We denote $\omega_0=0,\;\omega_3=\omega_1+\omega_2$,
then $e_j=\wp(\omega_j/2)$, $1\leq j\leq 3$. The resulting differential equation
is
called the Heun equation
in the elliptic form:
\begin{equation}\label{heun}
w''=\left(\sum_{j=0}^3 k_j\wp(z-\omega_j/2)+\lambda\right) w, \quad
k_j=\alpha_j^2-1/4.
\end{equation}
The two parameters are now the modulus of the torus $\tau=\omega_2/\omega_1$
and the
accessory parameter $\lambda$. The relation between $q$ and $\lambda$ is
affine: $\lambda=aq+b$ with $a$ and $b$ depending on $\tau$.
The precise form of $a$ and $b$ is irrelevant here,
but we will use the fact that they are both real when $a_j$ and $\alpha_j$
are real.
The details of calculation reducing (\ref{heun0}) to (\ref{heun})
are given in \cite{Tak}.

We will use both forms of the Heun equation. Equation (\ref{heun})
is considered on the torus (not on the plane). In particular
the monodromy includes the translations of solutions
by the elements of the lattice.

The ratio of two solutions $f=w_1/w_2$
is a developing map of a metric in question if and only if the projective
monodromy group
of the Heun equation is conjugate to a subgroup of $PSU(2)\simeq SO(3)$.
In this case
we say that the monodromy is {\em unitarizable}. The exponents
at the singularity $\omega_j/2$ are $\rho_j^\pm=1/2\pm\sqrt{1/4+k_j}$,
so the angle at this singularity is $4\pi\sqrt{1/4+k_j}$,
which is $4\pi\alpha_j$, twice the angle of the original metric on
the sphere.

The problem is for given $\omega_1,\omega_2$ and $k_j>-1/4,\; 0\leq j\leq 3$,
to find the values of $\lambda$ for which the monodromy is unitarizable.

\section{Statement and discussion of results}

The main result of this paper is
\begin{theorem}
For every $\tau$ and $k_j>-1/4$, $0\leq j\leq 3$, the set 
$$U:=\{\lambda\in\C: \mbox{equation (\ref{heun}) has unitarizable projective monodromy}\}$$
is finite. 
\end{theorem}

The correspondence between the metrics and Heun's equations is not
one-to-one:
different pairs of linearly independent solutions of the same equation
may correspond to different metrics. This can only happen when the
projective monodromy is {\em co-axial} that is isomorphic to
a subgroup of the unit circle. We say that a metric is
co-axial if the monodromy
of its developing map is co-axial. Co-axial metrics come in
continuous families consisting of {\em equivalent metrics}:
two metrics with developing maps $f_1,f_2$ are called equivalent
if $f_1=\phi\circ f_2$ for some linear-fractional transformation $\phi$.
When the projective
monodromy is trivial there is a real $3$-parametric 
family of metrics, and when the projective monodromy
is a non-trivial subgroup of
the circle there is a real $1$-parametric family of metrics,
see for example \cite{CWWX}.

Co-axial metrics on the sphere have been completely described in \cite{E2};
in particular, when $n\geq 3$ some angles of
a co-axial metric must be integer multiples of $2\pi$.
So Theorem 1 has the following
\begin{corollary}
For every four points $a_0,\ldots,a_3$ on the Riemann 
sphere and  every $\alpha_j\in\R_+\backslash\Z$, there exist at most
finitely many
metrics of curvature $1$ and conic singularities at $a_j$ with angles 
$2\pi\alpha_j$.
\end{corollary}

Whether a metric on the sphere is co-axial or not is completely determined
by the angles \cite[Theorem A]{MP1}. The only co-axial metrics on surfaces
of genus $g\geq 1$ are metrics on tori with all angles multiples of $2\pi$
\cite[Theorem 2.1]{MP2}. 

\begin{corollary} For every non-integer $\alpha$, and every torus, there
are at most finitely many metrics of curvature $1$ and one conic
singularity with angle $\alpha$ on this torus.
\end{corollary}

Indeed, the developing map of such a metric is a ratio of two
linearly independent solutions of (\ref{heun}) with $k_1=k_2=k_3=0$.

It follows from the results of Chen and Lin \cite{CL} that
the number of metrics on a torus with one singularity with
angle $\alpha$ is at least $[(\alpha-1)/2]+1$ unless $\alpha$
is an odd integer.

Unfortunately, our proof is non-constructive, and does not give
any explicit upper estimate for the number of metrics. 
The proof of Theorem 1 consists of two parts: first we prove that the set $U$
is bounded;
this part is based on the asymptotic analysis of equation (\ref{heun})
as $\lambda\to\infty$.

Compactness of the set of metrics with prescribed angles
in the given conformal class
has been proved in \cite{MP2} for metrics
on arbitrary compact Riemann surfaces with any number of singularities,
and for generic angles. For the case of $4$ singularities
on the sphere the condition on the angles in \cite{MP2}
is the following:
\vspace{.1in}

{\em None of the sums $\sum_{j=0}^3\pm\alpha_j$ is an even non-zero integer.}
\vspace{.1in}

This result also follows from \cite[Corollary 3]{Ba2} and
\cite[Corollary 2.5]{Ba3} with proofs based on a different method. 

The second part of our proof shows that the set of accessory parameters
defining unitarizable monodromy is discrete. A general theorem of Luo
\cite{L} implies that equations with unitarizable monodromy correspond
to a real analytic surface in the complex two dimensional space
of all Heun's equations with prescribed exponents. Assigning the position
of singularities means that we take the intersection of this real
analytic surface with a complex line. In general, 
such an intersection does not have to be discrete.

Finding the accessory parameters
corresponding to unitarizable monodromy requires solving
a system of equations of
the form
\begin{equation}\label{harm}
g_j(\lambda)=0,\quad j=1,2,
\end{equation}
where $\lambda$ is a complex variable and $g_j(\lambda)$ are real
harmonic functions,
and there is no general method
of proving that the set of solutions of (\ref{harm}) is discrete,
or to estimate the number of solutions from above. See \cite{LW}, \cite{BE}
where a very special case is solved.

To investigate equation (\ref{harm}) in our case,
we use a general
theorem of Stephenson \cite{St} which reduces the local question
about discreteness of the set of solutions of (\ref{harm}) to a question
about asymptotic behavior at infinity of entire functions 
(traces of the generators of monodromy),
and this question is
solved using the asymptotic behavior established in the first part
of the proof
and a theorem of Baker \cite{B} on compositions of entire functions.

It is the use of Stephenson's theorem that prevents our method
from working for $n>4$. In general, the set $U$ as in Theorem 1
is defined as a common zero set of
$2n-6$ real harmonic functions:
$$ g_j(\lambda)=0,\quad 1\leq j\leq 2n-6$$
of $n-3$ complex variables $\lambda=(\lambda_1,\ldots,\lambda_{n-3})$.
When $n=4$ we have two harmonic functions of one complex variable.
If the set of common zeros is not discrete, it must be unbounded.
This permits to use asymptotics as $\lambda\to\infty$ to obtain a contradiction.
But when $n>4$, the set of $2n-6$ real harmonic equations in
$n-3$ variables may have
a bounded non-discrete set of solutions,
so our argument does not work.

Nevertheless we
state the following
\begin{conjecture}
On any compact Riemann surface, the set of Fuchsian equations
with prescribed singularities and exponents, and with unitarizable
monodromy is finite.
\end{conjecture}
The problem considered here is somewhat similar to the accessory
parameter problem studied by Klein and Poincar\'e; see \cite{SG} for
a modern exposition of their work. They tried
to prove that there is a unique choice
of accessory parameters such that the ratio of solutions of
a Fuchsian differential equation is the inverse to the uniformizing map
of $S$ minus the punctures.
They did not succeed in proving the Uniformzation theorem with
this approach,
but the proof based on these ideas is completed in \cite{SG}.
It is interesting to notice that Poincar\'e did obtain a complete proof
for the case of the sphere with four punctures \cite{Po}.
This work of Poincar\'e
was continued by V. I. Smirnov \cite{S}, \cite{S1} whose argument is
used in Section 3 below. The main differences between our problem
and the problem of Klein and Poincar\'e
is that our problem may have more than one
solution, and that the monodromy group in our case is not discrete.

The author thanks Walter Bergweiler, Chang-Shou Lin,
Andrei Gab\-ri\-e\-lov, and Vitaly Tarasov for
useful discussions and comments.

\section{Asymptotics of traces of monodromy generators
and boundedness of $U$}

\def\tr{\mathrm{tr}\,}
Let $z_0$ be a point such that no singularities of equation (\ref{heun})
lie on the lines $$L_j=\{ z_0+t\omega_j: t\in\R\},\; j=1,2.$$
We restrict equation (\ref{heun}) on the line $L_j$ and obtain
an equation with periodic analytic potential with period $\omega_j$.
Let $T_j(\lambda)$ be the monodromy transformation corresponding
to the translation by $\omega_j$. Then the well-known result \cite{Mar},
\cite[Theorem 1]{Tka}
says that the trace $\tr T_j$ is an even entire function of $\sqrt{\lambda}$
with
the following asymptotic behavior:
\begin{equation}\label{trace}
\tr T_j(\lambda)=(2+O(1/\lambda))\cosh(\omega_j\sqrt{\lambda}),\quad\lambda\to
\infty,\quad|\arg(\omega^2\lambda)|<\pi-\epsilon,
\end{equation}
for every $\epsilon>0$.
These traces are usually called
Hill's discriminants or Lyapunov's functions in the literature on the
Sturm-Liouville equations.

If the monodromy is unitarizable, both traces must satisfy
$$\tr T_j(\lambda)\in[-2,2],$$
which is inconsistent with (\ref{trace}) for large $\lambda$
since the ratio $\omega_2/\omega_1$ is not real.
This proves that the set $U$ in Theorem 1 is bounded.

\section{The real case}

To prove discreteness of the set $U$ we first address the real
case: we assume that $(a_0,a_1,a_2)=(0,1,t)$, in (\ref{heun0})
and that $t$ and $q$ are real. Here we essentially follow
Smirnov \cite{S}, \cite{S1}, who investigated the
case of $SL(2,\R)$ monodromy and $\alpha_j<1$.

We only consider the case when none of the $\alpha_j$ in (\ref{heun0})
is an integer. For the case of at least one integer $\alpha_j$
our result was proved in \cite{ET}.

We assume without loss of generality that $t<0$. Let 
$w_{01},w_{02}$ be solutions of (\ref{heun0}) normalized 
at $0$ by
$$w_{0j}(z)=z^{\rho_{0j}}(1+g_{0j}(z)),$$
$\rho_{0j}$ are the exponents at $0$, $\rho_{00}=0,\; \rho_{01}=\alpha_0$,
and $g_{0j}$ are holomorphic and
vanish at $0$. Here and in what follows we use the principal
branches of powers. These two solutions $w_{0j}$ are 
real on $(0,1)$. We also consider two solutions 
$w_{11},w_{12}$ which are normalized
at $1$ and both are real on $(0,1)$:
$$w_{1j}(z)=(1-z)^{\rho_{1j}}(1+g_{1j}(z)),\quad z\in(1-\epsilon,1),$$
where $g_{1j}$ is analytic near $1$, $g_{1j}(1)=0$, and $\rho_{1j}$ are
the exponents at $1$, $\rho_{11}=\alpha_1,\; \rho_{12}=0$.

Then we have the connection matrix 
$$F=\left(\begin{array}{ll}f_{11}&f_{12}\\ f_{21}&f_{22}\end{array}\right),$$
such that
\def\w{\mathbf{w}}
$$\w_0=F\w_1,\quad\mbox{where}\quad\w_i=\left(\begin{array}{c}w_{i1}\\ w_{i2}
\end{array}\right),\quad i\in\{0,1\}.$$
The entries of $F$ are entire functions of the accessory parameter $q$,
real on the real
line, see \cite{S}, \cite{S1}. To obtain the projective monodromy,
we consider $f_i=w_{i1}/w_{i2},$
$i\in\{0,1\}$, which are related by a linear-fractional transformation
$f_0=L(f_1)$ which is represented by the matrix $F$.
Projective monodromy of $f_0$ at $0$ is an elliptic transformation
with fixed points $0$ and $\infty$. Since we assume that
none of the $\alpha_j$ is an integer, the local monodromy
at the singularities is not identical.
Monodromy of $f_0$ at $z=1$
is an elliptic transformation with fixed points $u_1=f_{11}/f_{21}$
and $u_2=f_{12}/f_{22}$. These points are real. Projective monodromies
at $0$ and $1$ are simultaneously unitarizable if and only if
the product of these fixed points is negative. Indeed,
an elliptic transformation is a rotation of the Riemann sphere if and
only if its fixed points $u_1,u_2$ are diametrally opposite,
that is 
\begin{equation}\label{do}
u_1\overline{u_2}=-1.
\end{equation}
In our case both $u_1,u_2$ are real so the bar can be dropped.
Choosing the fixed points of projective monodromy at $0$ to be $0,\infty$,
we still can multiply $f_0$  by a constant $\mu$. This will
result in multiplying both fixed points of the projective monodromy
at $1$ by $\mu$, so (\ref{do}) can be achieved for these fixed points
if and only if $u_1u_2<0$.

Similar considerations apply to the interval $(t,0)$. If we denote
the connection matrix on $(t,0)$ by
$$G=\left(\begin{array}{ll}g_{11}&g_{12}\\ g_{21}&g_{22}\end{array}\right),$$
then $g_{ij}$ are entire functions on $q$, real on the real line,
and the fixed points of the projective monodromy of $f_0$ at $t$
are $v_1=g_{11}/g_{21}$ and $v_2=g_{12}/g_{22}$.

So the condition of unitarizability is
\begin{equation}\label{me}
\frac{f_{11}f_{12}}{f_{21}f_{22}}=\frac{g_{11}g_{12}}{g_{21}g_{22}}<0.
\end{equation}
This includes the condition that two meromorphic functions of $q$
take equal values, so the set of $q$ satisfying this condition
is discrete, unless the equality in (\ref{me}) is
satisfied identically.

To show that the equation in (\ref{me}) cannot be satisfied
identically in $q$, one can use the asymptotics of solutions
for large $\lambda$,
but an easier way is to see this is directly
from the Heun equation in the form (\ref{heun0}), which in our case 
can be written as
$$w''+p(z)w'+q(z)w=0,$$
where
$$q(z)=\frac{Az-q}{z(z-1)(z-t)},\quad t<0.$$
When $q$ is large negative, $q(z)$ is large negative on $(0,1)$, so
solutions oscillate on $(0,1)$,
and functions $f_{ij}$ have infinitely many positive zeros,
while on the interval $(t,0)$ we have $q(z)>0$,
solutions do not oscillate, and functions $g_{ij}$ have no
large positive zeros.  Thus (\ref{me}) cannot
hold identically, and the set of real $q$ for which the monodromy
is unitarizable is discrete.

\section{Completion of the proof of Theorem 1}

To prove the second part of Theorem 1, discreteness of the set $U$,
we consider
two entire functions $\lambda\mapsto\tr T_j(\lambda)$ introduced in section 2.
They are the traces of the generators of the monodromy corresponding to
the periods $\omega_1,\omega_2$.

If for some $\lambda$
the monodromy is unitarizable, then $\tr T_j(\lambda)\in[-2,2]$,
so if the set of such $\lambda$ is not discrete,
there is a non-degenerate curve
$\gamma$ such that both $\tr T_j$ are real on $\gamma$.

Now we use the following
\vspace{.1in}

\noindent
{\bf Theorem of Stephenson} \cite[Thm. 13]{S}.
{\em Let $g_j,\; j=1,2$ be two entire functions which are both
real on a non-degenerate curve $\gamma$. Then
\begin{equation}\label{fac}
g_j=G_j\circ\phi,
\end{equation}
where $\phi,\; G_j$ are entire, $G_j$ are real on the real line,
and $\phi$ is real on $\gamma$.}
\vspace{.1in}

Recalling  the asymptotics (\ref{trace}) we obtain
\begin{equation}\label{ass2}
g_j(\lambda):=\tr T_j(\lambda)\sim 2\cosh(\sqrt{\omega_j^2\lambda}),
\quad\lambda\to\infty,\quad |\arg\omega^2\lambda|\leq\pi-\epsilon.
\end{equation}
These two functions have two {\em different} directions of maximal
growth and their zeros have arguments accumulating in the directions
opposite to the directions of maximal growth. To be more precise,
we say that an entire function $g$ has a single direction of
maximal growth $\theta$ if for every $\epsilon>0$ there exists $\delta>0$
such that for all $r>r_0$ we have
$$\max\{|g(re^{it})|:\theta+\epsilon\leq t\leq\theta+2\pi-\epsilon\}\leq
(1-\delta)\{\max|g(z)|:|z|= r\}.$$
Each of our functions $g_j$, $j=1,2$, has a single
direction of maximal growth $\theta_j=\arg(\omega_j^{-2})$,
and these directions are distinct because $\omega_1/\omega_2$
is not real. It follows that $G_j$ in (\ref{fac}) cannot be polynomials,
and $\phi$ cannot be a polynomial of degree greater than $1$.
Now we use 
\vspace{.1in}

\noindent
{\bf Theorem of Baker} \cite{B}.
{\em If an entire function $g$ of finite order
has a representation (\ref{fac}) with some entire transcendental
functions $G$ 
and all zeros of $g$ except finitely many lie in a sector of opening less
than $\pi$,
then $\phi$ must be a polynomial of degree~$1$.}
\vspace{.1in}

From this we conclude that the curve $\gamma$ in Stephenson's
theorem must be an interval of a straight line $\ell$,
and both $g_j$ are symmetric with
respect to this line, that is $g_j\circ s=\overline{g_j}$,
where $s$ is the reflection with respect to $\ell$.
Comparing this with asymptotics (\ref{ass2}) we conclude that
the directions of maximal growth of the $g_j$ must be collinear, and
since they are distinct, they must be opposite,
which gives $\omega_1=i\omega_2$. So our torus must be rectangular.
In terms of equation (\ref{heun0}), this means that
the singularities are real. This implies that our functions $g_j$
are real on the real line. 

Now we prove that $\ell$ is the real line. First $\ell$ cannot cross the
real line, because a function with two lines of symmetry will have at least
two directions
of maximal growth, while our functions have only one. Second,
it cannot be parallel to the real line, because in this case our functions
$g_j$ would be periodic which is incompatible with their asymptotics (\ref{ass2}).

This reduces the general case to the real case considered in the
previous section and completes the proof of Theorem 1.

\vspace{.2in}

{\em Purdue University, West Lafayette, IN 47907 USA

eremenko@math.purdue.edu}
\end{document}